\documentclass[12pt]{iopart}
\usepackage{psfrag}
\usepackage{subfigure}
\usepackage{graphicx}
\usepackage{iopams}
\def\strutdepth{\dp\strutbox}
\def\nw#1{\strut\vadjust{\kern-\strutdepth\vtop to0pt{\vss\hbox to\hsize
{\hskip\hsize\hskip5pt$\leftarrow$\hss\strut}}}{\em #1}}

\begin{document}

\title[Singularities]
{A catalogue of singularities}

\author{Jens Eggers$^*$ and Marco A. Fontelos$^\dagger$}

\address{
$^*$School of Mathematics,
University of Bristol, University Walk, \\
Bristol BS8 1TW, United Kingdom \\
$^\dagger$ Instituto de Matem\'aticas y F\'{\i}sica Fundamental,\\
Consejo Superior de Investigaciones Cient\'{\i}ficas \\
C/ Serrano 121, 28006 Madrid, Spain.
  }

\begin{abstract}
This paper is an attempt to classify finite-time singularities
of PDEs. Most of the problems considered
describe free-surface flows, which are easily observed experimentally.
We consider problems where the singularity occurs at a point, and where
typical scales of the solution shrink to zero as the singularity is
approached. Upon a similarity transformation, exact
self-similar behaviour is mapped to the fixed point of a
{\it infinite dimensional dynamical system} representing the original
dynamics. We show that the dynamics
close to the fixed point is a useful way classifying the structure of
the singularity.
Specifically, we consider various types of stable and unstable
fixed points, centre-manifold dynamics, limit cycles, and chaotic
dynamics.
\end{abstract}

\maketitle

\section{Introduction}
A non-linear partial differential equation (PDE),
starting from smooth initial data, will in general
not remain smooth for all times.
Consider for example the physical case shown in Fig.~\ref{lava},
which we will treat in section \ref{travel} below. Shown is
a snapshot of one viscous fluid dripping into another fluid,
close to the point where a drop of the inner fluid pinches off.
This process is driven by surface tension, which tries to minimise
the surface area between the two fluids. At a particular point
$x_0,t_0$ in space and time, the local radius $h(z,t)$ of the
fluid neck goes to zero. This point is a singularity of the
underlying equation of motion for the one-dimensional profile
$h(z,t)$. In particular, typical length scales near the pinch
point go to zero (the minimum radius, at the very least, but
also the solution's axial extend). This absence of a characteristic
scale near the singularity is the basic motivation to look for
self-similar solutions.

\begin{figure}[t]
\centering
\includegraphics[width=0.2\hsize]{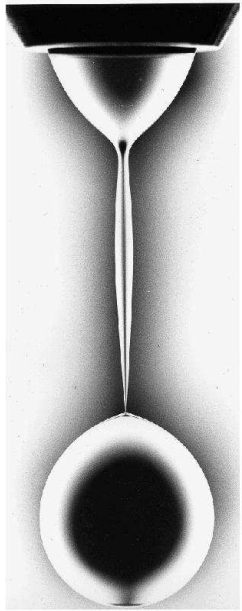}
\caption{\label{lava}
A drop of Glycerin dripping through PolyDimethylSiloxane
(PDMS) near snap-off \cite{CBEN99}. The nozzle diameter is $0.48$ cm.
   }
\end{figure}

A fascinating aspect of the study of singularities is that they
describe a great variety of phenomena which appear in the natural
sciences and beyond \cite{K97}. For example, such singular events
occur in free-surface flows \cite{E97}, turbulence and Euler dynamics
(singularities of vortex tubes \cite{M00,GMG98} and sheets \cite{CFMR05}),
elasticity \cite{AB03},
Bose-Einstein condensates \cite{BR02}, non-linear wave physics \cite{MGF03},
bacterial growth \cite{HV96,BCKSV99}, black-hole cosmology \cite{C93,MG03},
and financial markets \cite{Sor03}.

In this paper we consider equations
\begin{equation}
h_t=F[h],
\label{ge}
\end{equation}
where $F[h]$ represents some (nonlinear) differential or integral operator.
For simplicity, we will for the most part discuss the case of scalar
$h$, but make reference to the many important cases where the variable $h$
is a vector, or (\ref{ge}) is a system of equations.
Let us suppose that (\ref{ge}) forms a localised singularity at
$x_0,t_0$. If $t' = t_0-t$ and $x' = x-x_0$, we are looking for local
solutions of (\ref{ge}) which have the structure
\begin{equation}
h(x,t)= t'^{\alpha}\phi(x'/t'^{\beta}),
\label{ss}
\end{equation}
with appropriately chosen values of the exponents $\alpha,\beta$.

Giga and Kohn \cite{GK85,GK87} proposed to introduce self-similar
variables $\tau=-\ln(t')$ and $\xi=x'/t'^{\beta}$ to study the
asymptotics of blow up. Namely, putting
\begin{equation}
h(x,t)=t'^{\alpha}H(\xi,\tau),
\label{simg}
\end{equation}
(\ref{ge}) is turned into the ``dynamical system''
\begin{equation}
H_{\tau}=G[H]\equiv \alpha H+ t'^{1-\alpha}F[t'^{\alpha} H].
\label{ds}
\end{equation}
If (\ref{ss}) is indeed a solution of (\ref{ge}),
the right hand side of (\ref{ds}) is {\it independent} of
$\tau$, and self-similar solutions of the form (\ref{ss}) are
{\it fixed points} of (\ref{ds}). By studying the ``long-time''
($\tau\rightarrow\infty$) behaviour of solutions of (\ref{ds})
one can study the behaviour near blow-up. For the relation
of singular PDE problems to the renormalisation group,
developed in the context of critical phenomena, see
\cite{Goldenfeld93,BKL94}, for a more computational perspective,
see \cite{CDGK04}.

Solutions to the original PDE (\ref{ge}) for given initial data
can be viewed as orbits in some infinite dimensional phase phase,
for instance, $L^2$. If the fixed point is an attractor, blow-up
will be self-similar for some class of initial initial conditions.
However, other types of attractors ($\omega$-limit sets in
the notation which is customary in the context of partial
differential equations, see \cite{GV04} and references therein)
are frequently observed, furnishing
a very fruitful means of {\it classifying} singularities. In
this paper we discuss the following cases:

\begin{enumerate}
\item[(i)] {\it Stable fixed points}

Initial conditions in some neighbourhood of the fixed point
become attracted to it, so the solution converges exponentially
to the self-similar solution (\ref{ss}). Formally, two eigenvalues
of the linearisation around the fixed point are always positive,
but these unstable motions can be absorbed into a redefinition
of the origin of space and time; this is discussed in section
\ref{sub:first}. A sub-classification
into self-similarity of the first and of the second kind is due
to Barenblatt \cite{B96}. Self-similar solutions are
of the first kind if (\ref{ss}) only solves (\ref{ge}) for one
set of exponents $\alpha,\beta$; their values are fixed by either
dimensional analysis or symmetry, and are thus rational.
Solutions are of the second kind (in the sense of Barenblatt)
if solutions (\ref{ss}) exist for a continuous set of exponents
$\alpha,\beta$; the exponents are fixed by a non-linear eigenvalue
problem, and take irrational values in general.

\item[(ii)] {\it Unstable fixed points}

The linearisation around the fixed point possesses positive
eigenvalues, so it is never reached, except for non-generic
initial data. Often a stable fixed point is associated with an
infinite sequence of unstable fixed points.

\item[(iii)] {\it Travelling waves}

Solutions of (\ref{ge}) converge to
$h=t'^{\alpha}\phi(\xi+c\tau)$, which is a travelling wave solution
of (\ref{ds}) with propagation velocity $c$.

\item[(iv)] {\it Centre manifold}.

This is also known, when it leads to singularities that develop at
a faster rate than the selfsimilar scaling, as type-II self-similarity
\cite{AV97}; it arises if one
eigenvalue of the linearisation around the fixed point is zero,
and there is a non-linear dependence instead. This typically
leads to corrections involving logarithmic time $\tau$.
We discuss two different cases, involving quadratic and cubic
non-linearities.

\item[(v)] {\it Limit cycles}

This is also known as ``discrete self-similarity'' \cite{C93,MG07}.
Corresponding solutions have the form
$h={t^{\prime}}^{\alpha }\mathbf{\psi }\left[\xi,\tau\right]$
with $\psi$ being a periodic function of period $T$ in $\tau$.
Thus at the discrete sequence of times
$\tau_{n}=\tau_{0}+n\mathbf{T}$, which approaches the
singular time for $n\rightarrow \infty $, the solution looks
like a simple self-similar one.

\item[(vi)] {\it Strange attractors}

In principle, more complex behaviour is possible, where the orbits
of the dynamical system lie on a strange attractor. At the moment,
we are not aware of an equation exhibiting such behaviour which
would correspond to any physical phenomenon. However, this may
simply be due to the fact that the corresponding singular behaviour
is more difficult to detect numerically, and to grasp analytically.
To demonstrate that such behaviour is at least possible, we show that
{\it any} finite dimensional dynamical system may be ``embedded'' in the
singular dynamics. As an explicit example, we show that the phase-space
trajectory may lie on the Lorenz attractor.

\item[(vii)] {\it Exotic objects}

There might be other types of behaviour that have no analogue in
finite-dimensional dynamical systems. In particular, blow-up
may occur at several points $(x_0,t_0$ at the same time, in which
case the description (\ref{ds}) is not so useful.

\end{enumerate}

This paper's aim is to assemble the body of knowledge on
singularities of equations of the type (\ref{ge}) that is
available in both the mathematical and the applied community,
and to categorise it according to the types given above.
In addition to rigorous results we pay particular attention to
various phenomenological aspects of singularities which are often
crucial for their appearance in an experiment or a numerical simulation.
For example, what are the implications of the type of singularity
for the approach of the PDE solution onto the self-similar form
(\ref{ss})? In most cases, we rely on known examples from the
mathematical physics literature. To find an explicit example for
limit cycle behaviour as well as chaotic dynamics, we propose a
new set of model equations, inspired by a problem in general
relativity.

\section{Stable and unstable fixed points}\label{fixed}
\subsection{Self-similarity of the first kind}\label{sub:first}
\begin{figure}[t]
\centering
\includegraphics[width=0.4\hsize]{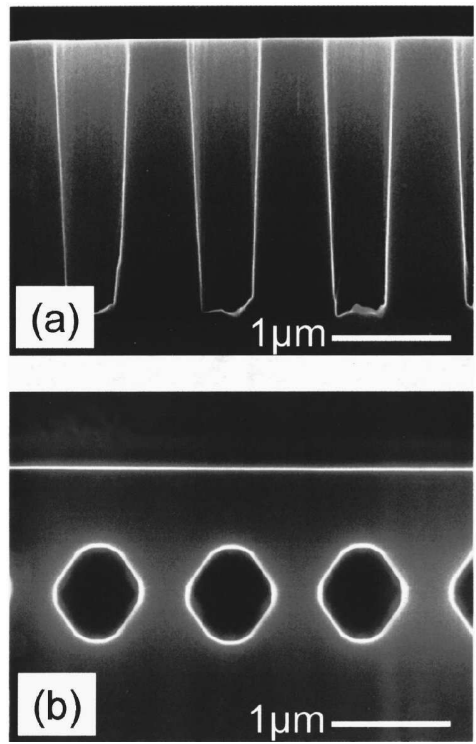}
\caption{\label{trough}
Scanning electron microscopy images illustrating the pinch-off of a row of 
rectangular
troughs in silicone (top) \cite{MSTT00}. The bottom picture
shows the same sample after 10 minutes of annealing at $1100^{\circ}$C.
The troughs have pinched off to form a row of almost spherical
voids. The dynamics is driven by surface diffusion.
   }
\end{figure}

Our example, exhibiting self-similarity of the first kind (in the sense
of Barenblatt) \cite{B96}, is that of a solid surface
evolving under the action of surface diffusion. Namely,
atoms migrate along the surface driven by gradients of
chemical potential, see Fig.\ref{trough}. The resulting equations in the
axisymmetric case, where the free surface is described
by the local neck radius $h(x,t)$, are \cite{NM65}:
\begin{equation}
h_t=\frac{1}{h}\left[\frac{h}{(1+h_x^2)^{1/2}}\kappa_x\right]_x,
\label{sd}
\end{equation}
where
\begin{equation}
\kappa = \frac{1}{h(1+h_x^2)^{1/2}} - \frac{h_{xx}}{(1+h_x^2)^{3/2}}
\label{mean}
\end{equation}
is the mean curvature. In (\ref{sd}),(\ref{mean}), all lengths
have been made dimensionless using an outer length scale $R$ (such
as the initial neck radius), and the time scale $R^4/D_4$, where
$D_4$ is a forth-order diffusion constant.

\begin{figure}
\centering
\includegraphics[width=0.6\hsize]{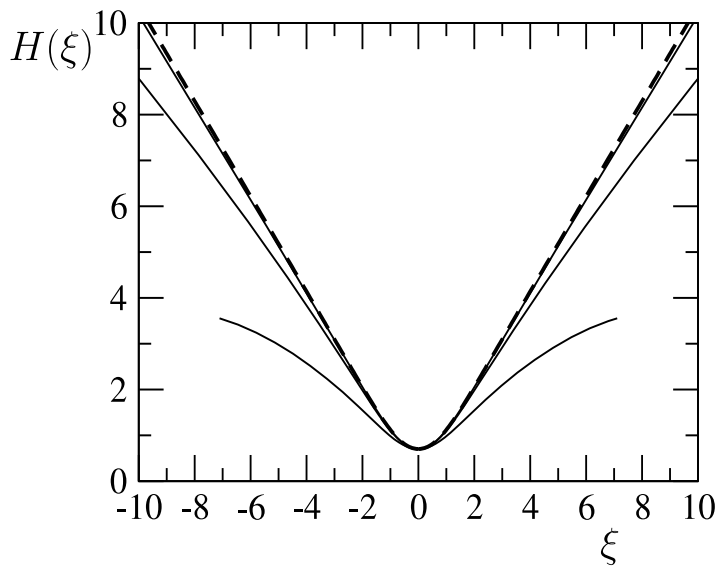}
\caption{\label{sdfig}
The approach to the self-similar profile for equation (\ref{sd}).
The dashed line is the stable similarity solution $\overline{H}(\xi)$
as found from (\ref{sdsim}). The full lines are rescaled profiles
found from the original dynamics (\ref{sd}) at
$h_{min} = 10^{-1},10^{-2}$, and $h_{min}=10^{-3}$, respectively.
As the singularity is approached, they converge rapidly onto the
similarity solution (\ref{ssd}).
   }
\end{figure}

At a time $t' \ll 1$ away from breakup, dimensional analysis
implies that $\ell=t'^{1/4}$ is a local length scale.
This suggests the similarity form
\begin{equation}
h(x,t)= t'^{1/4}\overline{H}(x'/t'^{1/4}),
\label{ssd}
\end{equation}
implying $\alpha=\beta=1/4$. This is the classical situation for
self-similarity of the first kind, more examples are found 
in \cite{B96,Edrop05}. The similarity form
of the PDE becomes
\begin{equation}
-\frac{1}{4}(\overline{H}-\xi \overline{H}') = \frac{1}{\overline{H}}
\left[\frac{\overline{H}}{(1+\overline{H}'^2)^{1/2}}
\overline{\kappa}'\right]', \quad \xi=\frac{x'}{t'^{1/4}}
\label{sdsim}
\end{equation}
where the prime denotes differentiation with respect to $\xi$.

Solutions of (\ref{sdsim}) have been studied extensively in \cite{BBW98}.
To ensure matching to a time-independent outer solution, the leading order
time dependence must drop out from (\ref{ssd}), implying that
\begin{equation}
\overline{H}(\xi) \sim c|\xi| , \quad \xi \rightarrow \pm\infty.
\label{bound}
\end{equation}
It turns out that all similarity solutions are symmetric, so only
one constant $c$ needs to be determined. Exactly as in the closely
related problem of surface-tension-driven fluid pinch-off, the
requirement of a certain growth condition (\ref{bound}) is enough
to fix a unique solution of (\ref{sdsim}) \cite{E93}. The
value of $c$ comes out as part of the solution. Solutions of
(\ref{sdsim}) with the growth condition (\ref{bound}) form a
discretely infinite set \cite{BBW98},
again like the fluids problem \cite{BLS96}. The series of similarity
solutions is conveniently ordered by descending values of the minimum,
see table \ref{series}.

\begin{table}
  \begin{center}
    \leavevmode
\begin{tabular}{cccccccccccccccccc}
i & $\overline{H}_i(0)$ & $c_i$ & \\ \hline
0 & 0.701595 & 1.03714 \\
1 & 0.636461 & 0.29866 \\
2 & 0.456842 & 0.18384 \\
3 & 0.404477 & 0.13489 \\
4 & 0.355884 & 0.10730 \\
5 & 0.326889 & 0.08942 \\
\end{tabular}
  \end{center}
\caption{A series of similarity solutions of (\ref{sdsim}) as given
in \cite{BBW98}. The higher-order solutions become successively thinner
and flatter. }
\label{series}
\end{table}

Next we turn to the dynamical system that describes the dynamics
away from the fixed point, by putting
\begin{equation}
h(x,t)=t'^{1/4} H(\xi,\tau),
\label{ssmean}
\end{equation}
where $\tau = -\ln(t')$. The similarity form of (\ref{sd}) becomes
\begin{equation}
H_{\tau}=\frac{1}{4}(H-\xi H_{\xi})+\frac{1}{H}
\left[\frac{H}{(1+H_{\xi}^2)^{1/2}}\kappa_{\xi}\right]_{\xi},
\label{mean_dyn}
\end{equation}
which reduces to (\ref{sdsim}) if the left hand side is set to
zero. To assure matching of (\ref{mean_dyn})
to the outer solution, we require the boundary condition
\begin{equation}
H_{\tau} - (H-\xi H_{\xi})/4\rightarrow 0 \quad \mbox{for}
\quad |\xi| \rightarrow \infty .
\label{bcsd}
\end{equation}

Next we linearise (\ref{mean_dyn}) around $H = \overline{H}(\xi)$,
by writing $H = \overline{H}(\xi) + \epsilon P(\xi,\tau)$, which gives
\begin{equation}
P_{\tau}={\cal L}(\overline{H}) P .
\label{sd_lin}
\end{equation}
Since $\overline{H}$ satisfies (\ref{bcsd}), $P(\xi,\tau)$ must do the
same. In particular, this means that if
\begin{equation}
{\cal L}(\overline{H}) P_i = \nu_i P_i,
\label{mean_lyn}
\end{equation}
i.e. if $\nu_i$ is an eigenvalue of the linear operator,
the corresponding eigenfunction must grow like
\begin{equation}
P_i(\xi) \propto \xi^{1-4\nu_i}.
\label{eigen_grow}
\end{equation}

If the similarity solution $\overline{H}(\xi)$ is to be stable,
the eigenvalues of ${\cal L}(\overline{H})$ must be negative.
However, there are always two positive eigenvalues, which are
related to the invariance of the equation of motion (\ref{sd})
under translations in space and time. Namely, for any $\epsilon$,
the translated similarity solution
\begin{equation}
h_{\epsilon}(x,t) = \overline{H}(\frac{x'+\epsilon}{\ell})
\label{trans}
\end{equation}
is an equally good self-similar solution of (\ref{sd}), and
thus of (\ref{mean_dyn}).
In particular, we can expand (\ref{trans}) to lowest order in
$\epsilon$, and find that
\begin{equation}
H_{\epsilon}(\xi,\tau) = \overline{H}(\xi) +
\epsilon e^{\beta\tau} \overline{H}'(\xi) + O(\epsilon^2)
\label{trans_sim}
\end{equation}
is a solution of (\ref{sd_lin}).

Thus, since ${\cal L}\overline{H}=0$,
\begin{equation}
\epsilon e^{\beta\tau} \beta \overline{H}'  =
\frac{\partial H_{\epsilon}(\xi,\tau)}{\partial\tau} =
{\cal L}H_{\epsilon} = \epsilon e^{\beta\tau} {\cal L} \overline{H}' .
\label{f_exp}
\end{equation}
But this means that $\nu_x = \beta\equiv 1/4$ is an eigenvalue of
${\cal L}$ with eigenfunction $\overline{H}'(\xi)$. Similarly, considering
the transformation $t\rightarrow t + \epsilon$, one finds a second
positive eigenvalue $\nu_t = 1$, with eigenfunction
$\xi\overline{H}'$. To reiterate, the physical meaning of these
eigenvalues is that upon perturbing the similarity solution, the singularity
time as well as the position of the singularity will change. Thus if
the coordinate system is not adjusted accordingly, it looks as if the
solution would flow away from the fixed point. If, on the other hand,
the solution is represented relative to the perturbed values of
$x_0$ and $t_0$, the dynamics will converge onto a stable similarity
solution.

The eigenvalues of the solutions $\overline{H}_i$ have been found
numerically in \cite{BBW98}. The result is that the linearisation
around the ``ground state'' solution $\overline{H}_0$ only has
negative eigenvalues (apart from the two trivial ones),
while {\it all} the other solutions have at least one other
positive eigenvalue. This means that $\overline{H}_0$ is the only
similarity solution that can be observed, all other solutions
are unstable. Close to the fixed point, the approach to $\overline{H}_0$
will be dominated by the largest negative eigenvalue $\nu_1$:
\begin{equation}
h(x,t)= t'^{1/4} \left[\overline{H}(\xi) +
\epsilon t'^{\nu_1} P_1(\xi)\right].
\label{approach}
\end{equation}
For large arguments, the point $\xi_{cr}$ where the correction becomes
comparable to the similarity solution is
$\xi \sim \epsilon t'^{\nu_1} \xi^{1-4\nu_1}$, and thus
$\xi_{cr} \sim t'^{1/4}$. This means that the region of validity
of $\overline{H}(\xi)$ {\it expands} in similarity variables, and
is constant in real space.
This rapid convergence is reflected by the numerical results reported in
Fig. \ref{sdfig}. More formally, one can say that for
any $\epsilon$ there is a $\delta$ such that
\begin{equation}
\left|h(x,t) - t'^{1/4}\phi(\xi)\right| \le \epsilon
\label{conv1}
\end{equation}
if $|x'|\le \delta$ {\it uniformly} as
$t' \rightarrow 0$.

\subsection{Self-similarity of the second kind}
In the example of the previous subsection, the exponents
can be determined by dimensional analysis, and therefore
assume rational values. As Barenblatt \cite{B96} points out,
there are problems where the scaling behaviour depends on
external parameters, set for example by the initial conditions.
In that case, the scaling exponent can assume any value.
Often, this value is fixed by some intrinsic property
of the equation, resulting in an irrational answer. We will call
this situation self-similarity of the second kind (in the sense of
Barenblatt).
A particularly simple example of this kind of singularity
is the pinch-off of a very viscous thread of liquid \cite{P95,E97},
which we present now. Another recent example is
the pinch-off of a two-dimensional inviscid sheet \cite{BT07}.

For simplicity, we confine ourselves to the case of a slender
viscous filament without inertia, for which the equation becomes:
\begin{equation}
h_{t}(s,t)=\frac{1}{6}\left(1+\frac{C(t)}{h(s,t)}\right).
\label{visc}
\end{equation}
The typical velocity scale $\gamma/\eta$, where $\gamma$ is the
surface tension and $\eta$ is the viscosity, has been absorbed
into the time variable. The particularly simple form of (\ref{visc})
has been achieved by writing the thread radius in
{\it Lagrangian variables}, i.e. as function of a particle
label $s$. This means the particle is at position $z(s,t)$
at time $t$, and $z_t(s,t)$ is the velocity at time $t$.
The time-dependent constant of integration $C(t)$ must be determined
from the constraint that
$u(s,t)\equiv z_{s }=1/h^{2}(s ,t)$.

Note that the self-similar form (\ref{ss}) is a solution of (\ref{visc})
for $\alpha = 1$, and {\it any} value of $\beta$. Thus contrary to the
example described in the previous section, the exponents are not
completely determined from dimensional analysis or from balancing
powers of $t'$ in the equation of motion. This is a typical
situation in which self-similarity of the second kind is observed
\cite{B86}. Instead, the unknown exponent is determined from a
non-linear eigenvalue equation, and takes an irrational value.

\begin{figure}[t]
\centering
\includegraphics[width=0.4\hsize]{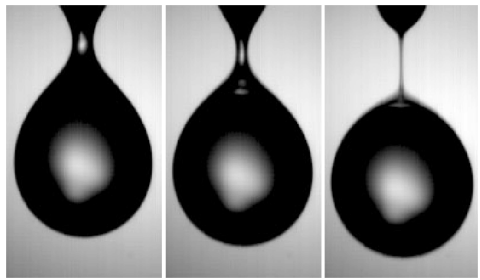}
\caption{\label{rrr}
A drop of viscous fluid falling from a pipette
1 mm in diameter \cite{RRR03}. Note the long neck.
   }
\end{figure}

Since $\alpha =1$ we introduce
\begin{equation}
u=t'^{-2}f\left( \xi \right),\quad \mbox{with}\quad
\xi =s/t'^{\beta}
\label{u}
\end{equation}
and
\begin{equation}
C(t)=K t' \ .
\label{paja21}
\end{equation}
Hence
\begin{equation}
\frac{1}{\sqrt{f}}+3\left( \frac{2}{f}+\frac{\beta \xi
f^{\prime }}{f^{2}}\right) =K  ,
\label{d1}
\end{equation}
where $K$ is an arbitrary constant. Imposing symmetry and
regularity of $f$, we introduce an expansion of $f(\xi )$ of
the form
\[
f(\xi )=R_{0}^{-2}+C\xi ^{2n}+O(\xi ^{4n})\ ,\ n=1,2,...
\]
into (\ref{d1}) to obtain the condition
\begin{equation}
R_{0}=\frac{1}{12(n\beta -1)}
\label{min}
\end{equation}
and define $\overline{\beta }=n\beta .$ Equation (\ref{d1}) can easily be
integrated in terms of $\ln \xi $ and $y=\sqrt{f}$:
\[
\int \frac{dy}{\left( \left( 1+6R_{0}\right)
y^{3}-y^{2}-6R_{0}y\right) }=\frac{1}{6R_{0}\beta }\ln \xi +
\widetilde{C}=\frac{1}{6R_{0}\overline{\beta }}
\ln \xi ^{n}+\widetilde{C},
\]
with $\widetilde{C}$ an arbitrary constant. Computing the
integral above we obtain
\begin{equation}
y^{-\beta }\left( \left( 2\overline{\beta}-1\right) y+
1\right) ^{\overline{\beta }-\frac{1}{2}}\left( 1-y\right)
^{\frac{1}{2}}=C\xi ^{n},
\label{d3.1}
\end{equation}
where we can fix, without loss of generality, $C=1$. In terms of
the profile $h(s,t)$ this corresponds to fixing the scale of the
spatial variable $s$. Solutions are undetermined up to such a scale
factor, as is clear from the invariance of (\ref{visc}) under a change
in spatial scale. As a result, the axial scale is fixed by the
initial conditions.

The value of the velocity $C_{n}$ at infinity is therefore given by
\[
C_{n}=\int_{0}^{\infty }z_{s t}d s=
\int_{0}^{\infty }u_{t}d s=\frac{1}{3}\int_{0}
^{\infty }\left( \left( \frac{1}{24}\frac{2\overline{\beta }-1}
{(\overline{\beta }-1)^{2}}\right) f^{2}-f^{\frac{3}{2}}\right) d\xi.
\]
From the condition that the velocity at infinity must vanish
we thus obtain $C_n = 0$, which will be the equation that
determines the exponent $\beta$.
Taking the derivative of (\ref{d3.1}) we obtain
\[
n\xi ^{n-1}\frac{d\xi }{dy}=\frac{d}{dy}
\left( y^{-\overline{\beta}}\left(\left( 2\overline{\beta }-1\right) y+
1\right) ^{\overline{\beta }-\frac{1}{2}}\left( 1-y\right)
^{\frac{1}{2}}\right) =
\]
\[
=-y^{-\overline{\beta }-1}
\left( 2y\overline{\beta }-y+1\right) ^{\overline{\beta }-
\frac{3}{2}}\frac{\overline{\beta }}{\sqrt{\left(1-y\right) }}
\]
and hence
\begin{eqnarray}
&& K_{n}(\beta )\equiv \frac{3}{(12(\overline{\beta}-1))^{3}}
C_{n}=\frac{\overline{\beta }}n\int_0^1
\left( \left( \frac 12\frac{2\overline{\beta }-1}{\overline{\beta }-1}\right)
y^4-y^3\right) \nonumber \\
&& \left( y^{-\frac{n+\overline{\beta }}n}
\left( \left( 2\overline{\beta }-1\right)y+1\right)^{-\frac 12
\frac{2n-2\overline{\beta }+1}n}\left( 1-y\right) ^{-\frac 12%
\frac{2n-1}n}\right) dy  .
\label{minne6}
\end{eqnarray}
The function $K_{n}(\beta )$ may be written explicitly as
\begin{eqnarray}
&& K_{n}(\beta ) =\beta \frac{\Gamma \left( 4-\beta \right)
\Gamma\left( \frac{1}{2n}\right) }{\Gamma
\left( 4-\beta+\frac{1}{2n}\right) }
\frac{2n\beta -1}{2n\beta -2} \nonumber \\
&& F\left( \frac{2n+1}{2n}
-\beta ,4-\beta ;4-\beta +\frac{1}{2n};1-2n\beta \right) - \nonumber \\
&&\ \ -\beta \frac{\Gamma \left( 3-\beta \right) \Gamma
\left( \frac{1}{2n}\right) }{\Gamma \left( 3-\beta +\frac{1}{2n}\right) }\;
F\left( \frac{2n+1}{2n}-\beta ,3-\beta ;3-\beta +
\frac{1}{2n};1-2n\beta \right)
\label{kn}
\end{eqnarray}
with roots given by table \ref{exponents}. If one converts the Lagrangian
variables back to the original spatial variables, one obtains
\begin{equation}
h(x,t) = t'\phi_{St}\left(x'/t'^{\beta-2}\right) .
\label{orig}
\end{equation}
Thus for $t'\rightarrow 0$ the typical radial scale $t'$ of the
generic $n=1$ solution rapidly becomes smaller than the axial scale
$t'^{0.175}$ (cf. table \ref{exponents}). This explains the long necks
seen in Fig.~\ref{rrr}.

\begin{table}
  \begin{center}
    \leavevmode
\begin{tabular}{ccccccccc}
n & $\beta$ & \\ \hline
1 & 2.1748 & \\
2 & 2.0454 & \\
3 & 2.0194 & \\
4 & 2.0105 & \\
5 & 2.0065 & \\
10 & 2.0014 & \\
\end{tabular}
  \end{center}
\caption{A list of exponents, found from $K_n(\beta) = 0$,
with $K_n$ given by (\ref{kn}).
The number $2n$ gives the smallest non-vanishing power in
a series expansion of the corresponding similarity solution
around the origin. Only the solution with $n=1$ is stable.
The minimum radius is found from (\ref{min}).
 }
\label{exponents}
\end{table}
For generic initial data
$u_{0}(s)=h_{0}^{-2}(s)=B_{0}+B_{1}s^2+B_{2}s^4+...+B_{j}s^{2j}+...$
one expects that $B_{1}\neq 0$, so that the self-similar solution with $n=1$
will develop. Only
if $B_{i}=0$ for $i=1,2,...,n-1$ and $B_{n}\neq 0$ the $n$-th self-similar
solution will be the asymptotic description of the solution. For this 
reason, only the $n=1$ solution is stable, since a generic perturbation 
of the initial
data with $B_{i}=0$ for $i=1,2,...,n-1$ will, in general, make $B_{1}\neq 0$.

\section{Travelling wave}
\label{travel}
The pinching of a liquid thread in the presence of an
external fluid is described by the Stokes equation \cite{LS98}. For
simplicity, we consider the case that the viscosity $\eta$ of the
fluid in the drop and that of the external fluid are the same.
An experimental photograph of this situation is shown in
Fig.~\ref{lava}.
To further simplify the problem, we make the assumption
(the full problem is completely analogous) that the fluid
thread is slender. Then the equations given in \cite{CBEN99} simplify to
\begin{equation}
\label{interface}
h_t = -v_z h/2 - v h_z,
\end{equation}
where
\begin{equation}
\label{velocity}
v = \frac{1}{4} \int_{z_-}^{z_+}
\frac{h_{z'}(z')}{(h^2(z')+(z-z')^2)^{1/2}} dz'.
\end{equation}
Here we have written the velocity in units of the capillary
speed $v_{\eta} = \gamma/\eta$. The limits of integration
$z_-$ and $z_+$ are for example the positions of the plates
which hold a liquid bridge \cite{P43}.

Dimensionally, one would once more expect a local solution of the
form
\begin{equation}\label{56}
    h(z,t)= t'H_{\rm{out}}\left(\frac{z'}{t'}\right),
\nonumber
\end{equation}
and $H_{out}(\xi)$ has to be a linear function at infinity to match
to a time-independent  outer solution. In similarity variables,
\eref{velocity} has the form
\begin{equation}\label{58}
    V_{\rm{out}}(\xi)=\frac
    14\int^{z_{b/t'}}_
      {-z_b/t'}\frac{H'_{\rm{out}}(\xi')}
    {\sqrt{H^2_{\rm out}+(\xi-\xi')^2}}d\xi' .
\end{equation}
We have chosen $z_b$ as a real-space variable close to the pinch-point,
such that the similarity description is valid in $[-z_b,z_b]$.
But if $H_{\rm{out}}$ is linear, the integral in (\ref{58}) diverges,
which means that a simple ``fixed point'' solution (\ref{56}) is
impossible.

However, the integral can be made convergent by introducing
a shift in the similarity variable $\xi$:
\begin{equation}\label{60}
    h=t' H_{\rm{out}}(\xi-b\tau),
\end{equation}
with $\tau=-\ln(t')$ as usual. This means
in similarity variables the solution is a travelling wave.
With this modification, the mass balance (\ref{interface}) becomes
\begin{equation}\label{61}
    -H_{\rm{out}}+H'_{\rm{out}}(\xi+V_{\rm{out}}+b\tau)
 =H_{out}V_{out}'/2. \nonumber
\end{equation}
Now we choose $b$ such that the logarithmic singularity cancels, namely
we demand that
\begin{equation}\label{62}
 \frac14\int^{z_b/t'}_{-z_b/t'}
\frac{H_{out}'(\xi')}{\sqrt{H_{out}^2+(\xi-\xi')^2}}d\xi'-
  b\ln(t') \nonumber
\end{equation}
finite for $t'\to 0$. This is achieved by putting
\begin{equation}\label{63}
  b=\frac14\left[-\frac{H_+}{\sqrt{H_+^2+1}}+
\frac{H_-}{\sqrt{H_-^2+1}}\right]. \nonumber
\end{equation}

Thus defining
\begin{equation}\label{64}
   V_{\rm{fin}}(\xi)=\lim_{\Lambda\to \infty}\frac 14
   \int^\Lambda_{-\Lambda}\frac{H'_{\rm{out}}}
{\sqrt{H^2_{\rm{out}}+(\xi-\xi')^2}}d\xi'+b\ln \Lambda
\end{equation}
the similarity equation
\begin{equation}\label{65}
  -H_{\rm{out}}+H'_{\rm{out}}(\xi+V_{\rm{fin}}-\xi_0)=
H_{\rm{out}}V'_{\rm{fin}}/2
\end{equation}
is finite, and $\xi_0$ is an arbitrary constant. It remains as
an arbitrary axial shift in the similarity solution.

The numerical solution of the integro-differential equation
(\ref{65}) gives
\begin{equation}\label{66}
    h_{\min}=a_{\rm{out}}v_\eta t', \quad \mbox{where} \quad
    a_{\rm{out}}=0.0335.
\end{equation}
The slope of the solution away from the pinch-point are
given by
\begin{equation}\label{slope}
 H_+=4.81 \quad \mbox{and} \quad  H_-=-0.105,
\end{equation}
which means the solution is very asymmetric, as confirmed
directly from Fig.~\ref{lava}.

\section{Centre manifold}
In section \ref{fixed} we described the generic situation
that the behaviour of a similarity solution is determined
by the linearisation around it. In the case of a stable fixed
point, convergence is fast, and the observed behaviour is essentially
that of the fixed point. In this section, we describe two
different cases where the largest eigenvalue vanishes, so higher-order
non-linear terms have to be taken into account. The approach to the
fixed point is now much slower, and much more of the observed behaviour
is determined by the approach to the fixed point. Depending on the
type of non-linearity, there is great freedom of possible behaviours,
of which we discuss two.

\subsection{Quadratic non-linearity}
\label{quad}
\begin{figure}
\centering
\centerline{\includegraphics[width=0.4\linewidth]{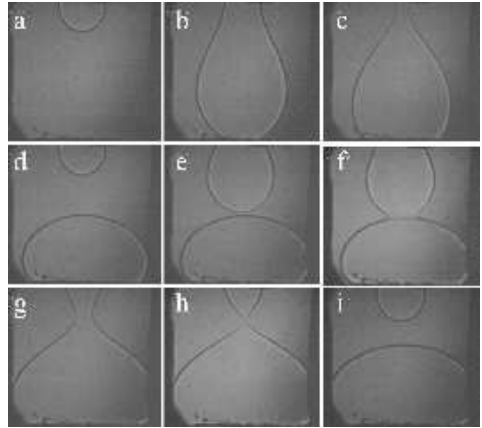}}
\vspace{0.5 cm}
\caption{Nine images (of width 3.5 mm) showing how
a $^3$He crystal ``flows'' down from the upper part of a cryogenic
cell into its lower part \cite{Ishiguro}. The recording takes a
few minutes, the temperature is 0.32 K.
11 mK. The crystal first ``drips'' down, so that a crystalline
``drop'' forms at the bottom (a to c); then a second drop appears
(d) and comes into contact with the first one (e); coalescence
is observed (f) and subsequently breakup occurs (h). }
\label{f2}
\end{figure}

Axisymmetric motion by mean curvature in three spatial
dimensions is described by the equation
\begin{equation}
h_t= \left(\frac{h_{xx}}{1+h_x^2} - \frac{1}{h}\right),
\label{mc}
\end{equation}
where $h(x,t)$ is the radius of the moving free surface.
A very good physical realization of (\ref{mc}) is the melting
and freezing of a $^3$He crystal, driven by surface tension \cite{IGRBE07},
see Fig.~\ref{f2}. As before, the time scale $t$ has been
chosen such that the diffusion constant, which sets the rate
of motion, is normalised to one.
A possible boundary condition for the problem is
that $h(0,t) = h(L,t) = R$, where $R$ is some prescribed
radius. For certain initial conditions $h(x,0)\equiv h_0(x)$
the interface will become singular at some time $t_0$,
at which $h(x_0,t_0)=0$ and the curvature blows up. The
moment of blow-up is shown in panel h of Fig.~\ref{f2},
for example.

Inserting the self-similar solution (\ref{ss}) into (\ref{mc}),
one finds a balance for $\alpha = \beta = 1/2$. The corresponding
similarity equation is
\begin{equation}
-\frac{\phi}{2}+\xi\frac{\phi'}{2} =
\left(\frac{\phi''}{1+\phi'^2} - \frac{1}{\phi}\right), \quad
\xi = \frac{x'}{t'^{1/2}}.
\label{mean_sim}
\end{equation}
One solution of (\ref{mean_sim}) is the constant
solution $\phi(\xi) = \sqrt{2}$.
Another potential solution is one that grows linearly at infinity,
to ensure matching onto a time-independent outer solution.
However, it can be shown that no solution to (\ref{mean_sim}),
which also grows linearly at infinity, exists \cite{AAG95,Huis93}.
Our analysis below follows the rigorous work in \cite{AV97},
demonstrating type-II self-similarity. In addition, we now show
how the description of the dynamical system can be carried out
to arbitrary order.

The solution relevant to our analysis is the constant solution, but which of
course does not match onto a time-independent outer solution. We thus
write the solution as
\begin{equation}
h(x,t)=t'^{1/2}\left[\sqrt{2}+g(\xi,\tau)\right],
\label{blmc}
\end{equation}
with $\tau = -\ln(t')$ as usual.
The equation for $g$ is then
\begin{equation}
g_{\tau}=g-\frac{\xi g'}{2}+\frac{g''}{1+g'^2} -
\frac{g^2}{2^{3/2}\sqrt{1+g'^2}},
\label{gequ}
\end{equation}
which we solve by expanding into eigenfunctions of the linear
part of the operator
\begin{equation}
{\cal L}g = g-\xi g'/2+g''.
\label{ling}
\end{equation}

It is easily confirmed that
\begin{equation}
{\cal L}H_{2i}(\xi/2) = (1-i)H_{2i}(\xi/2),
\label{eigen}
\end{equation}
where $H_n$ is the n-th Hermite polynomial \cite{AS68}:
\begin{equation}
H_n = (-1)^ne^{x^2}\frac{d^n}{dx^n}e^{-x^2}.
\label{Hermite}
\end{equation}
Thus the linear part of (\ref{gequ}) becomes
\begin{equation}
\frac{\partial a_i}{\partial\tau}=(1-i)a_i,
\label{linai}
\end{equation}
which means that all eigenvalues are negative except for the first, which
vanishes. To investigate the approach of the cylindrical solution,
one must therefore include nonlinear terms in the equation for $a_1$.

If we write
\begin{equation}
g(\xi,\tau) = \sum_{i=1}^{\infty} a_i(\tau)H_{2i}(\xi/2),
\label{gexp}
\end{equation}
the equation for $a_1$ becomes
\begin{equation}
\frac{d a_1}{d\tau} = -2^{3/2}a_1^2 + O(a_1a_j),
\label{a1}
\end{equation}
whose solution is
\begin{equation}
a_1 = 1/(2^{3/2}\tau).
\label{a1_r}
\end{equation}

Thus instead of the expected
exponential convergence onto the fixed point, the approach is only
algebraic. Since all other eigenvalues are negative, the $\tau$-dependence
of the $a_i$ is effectively determined by $a_1$. Namely, as we will see
below, $a_j=O(\tau^{-j})$, so corrections to (\ref{a1}) are of
higher order.

If one linearises around (\ref{a1_r}), putting $a_1 = a_1^{(0)} + \epsilon_1$,
one finds
\begin{equation}
\frac{d \epsilon_1}{d\tau} = -\frac{2}{\tau}\epsilon_1 + \mbox{other terms}.
\label{a1_pert}
\end{equation}
This means that the coefficient $A$ of $\epsilon_1 = A/\tau^2$ remains
undetermined, and a simple expansion of $a_i$ in powers of $\tau^{-1}$
yields an indeterminate system. Instead, at quadratic order, a term of
the form $\epsilon_1 = A\ln\tau/\tau^2$ is needed. Fortunately, this is
the only place in the system of nonlinear equations for $a_i$ where
such an indeterminacy occurs. Thus all logarithmic dependencies can be
traced, leading to the general ansatz
\begin{equation}
a_i^{(n)} = \frac{\delta_i}{\tau^i} + \sum_{k=i+1}^n
\sum_{l=0}^{k-i}\frac{(\ln\tau)^l}{\tau^k}\delta_{lki},
\label{ai_gen}
\end{equation}
where $\delta_i$ and $\delta_{lki}$ are coefficients to be determined.
The index $n$ is the order of the truncation.

Indeed, the coefficients can be found recursively by considering
terms of successively higher order in $\tau^{-1}$ in the first
equation:
\begin{eqnarray}
\label{ai_sys1}
&& \frac{d a_1}{d\tau} = -2^{3/2}a_1^2 - 24\sqrt{2}a_1a_2+22a_1^3 -\nonumber\\
&& 272\sqrt{2}a_1^4-191\sqrt{2}a_2^2+192a_1^2a_2 \\
\label{ai_sys2}
&& \frac{d a_2}{d\tau} = -a_2-\sqrt{2}/4a_1^2+6a_1^3-8\sqrt{2}a_1a_2
\end{eqnarray}
The next two orders will involve the next coefficient $a_3$.
From (\ref{ai_sys1}) and (\ref{ai_sys2}), one first finds
$\delta_{121}$ and $\delta_2$, by considering $O(\tau^{-3})$ and
$O(\tau^{-2})$, respectively. Then, at order $O(\tau^{-(n+1)})$
in the first equation, where $n = 3$, one finds all remaining
coefficients $\delta_{lki}$ in the expansion (\ref{ai_gen}) up
to $k=n$. At each order in $\tau^{-1}$, there is of course a
series expansion in $\ln\tau$ which determines all the coefficients.

We constructed a MAPLE program to compute all the coefficients up
to arbitrarily high order (10th, say). Up to third order in $\tau^{-1}$
the result is:
\begin{eqnarray}
&& a_1 = 1/4\,{\frac {\sqrt {2}}{\tau}}+{\frac {17}{16}}\,{\frac {\ln  \left(
\tau \right) \sqrt {2}}{{\tau}^{2}}}-{\frac {73}{16}}\,{\frac {\sqrt {
2}}{{\tau}^{3}}}+ \nonumber \\
&& {\frac {867}{128}}\,{\frac {\ln  \left( \tau \right)
\sqrt {2}}{{\tau}^{3}}}-{\frac {289}{128}}\,{\frac { \left( \ln
 \left( \tau \right)  \right) ^{2}\sqrt {2}}{{\tau}^{3}}} \\
&& a_2=-1/32\,{\frac {\sqrt {2}}{{\tau}^{2}}}+{\frac {5}{16}}\,{\frac {\sqrt
{2}}{{\tau}^{3}}}-{\frac {17}{64}}\,{\frac {\ln  \left( \tau \right)
\sqrt {2}}{{\tau}^{3}}},
\label{ai_res}
\end{eqnarray}
and thus $h(x,t)$ becomes
\begin{equation}
h(x,t)=t'^{1/2}\left[\sqrt{2}+a_1(\tau)\left(-2+\xi^2\right)
+ a_2(\tau)\left(12-12\xi^2+\xi^4\right)\right],
\label{h_res}
\end{equation}
from which one finds the minimum.
To second order, the result is
\begin{equation}
h_{min}=(2 t')^{1/2}\left[1 - \frac{1}{2\tau}
-\frac{3+17\ln\tau}{8\tau^2} \right].
\label{hmin_res}
\end{equation}

\begin{figure}
\centering
\includegraphics[width=0.5\hsize]{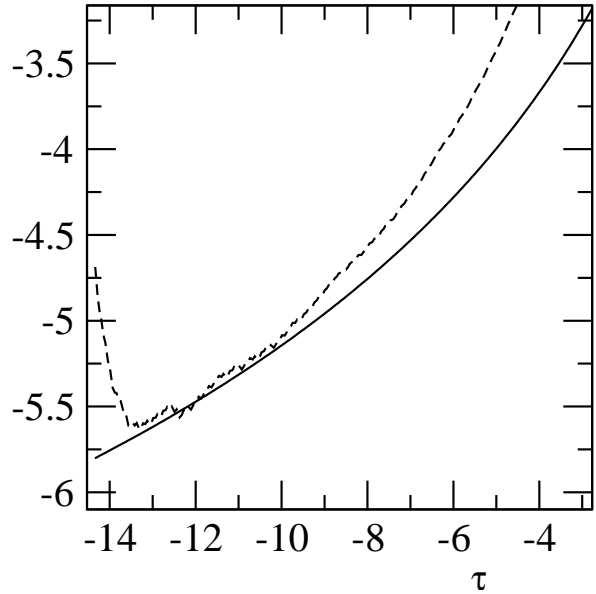}
\caption{\label{hminfig}
A plot of $\left[h_{min}/\sqrt{2t'}-1+1/(2\tau)\right]\tau^2$
(dashed line) and $\tau_0/2 - (3+17\ln(\tau+\tau_0)/8)$ (full
line) with $\tau_0 = 4.56$.
   }
\end{figure}
Two remarks are in order. First, the presence of logarithms implies
that there is some dependence on initial conditions built into
the description. The reason is that the argument inside the logarithm
needs to be non-dimensionalised using some ``external'' time scale.
More formally, any change in time scale $\tilde{t} = t/t_0$ leads to
an identical equation if also lengths are rescaled according to
$\tilde{h} = h/\sqrt{t_0}$. This leaves the prefactor in
(\ref{hmin_res}) invariant, but adds an arbitrary constant $\tau_0$
to $\tau$. This is illustrated by comparing to a numerical simulation
of the mean curvature equation (\ref{mc}) close to the point
of breakup, see Fig.~\ref{hminfig}. Namely, we subtract the analytical result
(\ref{hmin_res}) from the numerical solution $h_{min}/(2\sqrt{t'})$
and multiply by $\tau^2$. As seen in Fig.\ref{hminfig}, the remainder
is varying slowly over 12 decades in $t'$. If the constant $\tau_0$ is
adjusted, this small variation is seen to be consistent with the
logarithmic dependence predicted by (\ref{hmin_res}).

The second important point is that convergence in space is no longer
uniform as implied by (\ref{conv1}) for the case of self-similarity
of the first kind. Namely, to leading order the pinching solution
is a cylinder. For this to be a good approximation, one has to require
that the correction is small: $\xi^2/\tau \ll 1$. Thus corrections
become important beyond $\xi_{cr} \sim \tau$, which, in view of the
logarithmic growth of $\tau$, implies convergence in a constant region
{\it in similarity variables only}. As shown in \cite{IGRBE07}, the
slow convergence toward the self-similar behaviour has important
consequences for a comparison to experimental data.

\subsection{Cubic non-linearity}
The next example is that of bubble breakup \cite{EFLS07}, for
which a very different form of nonlinearity is observed. As shown
in \cite{EFLS07}, the equation for a slender cavity or bubble is
\begin{equation}
\int_{-L}^L\frac{\ddot{a}(\xi,t)d\xi}{\sqrt{(z-\xi)^2+a(z,t)}} =
\frac{\dot{a}^2}{2a} ,
\label{bernoulli}
\end{equation}
where $a(z,t)\equiv h^2(z,t)$. The integral runs over the
fluid domain.
If for the moment one disregards boundary conditions looks for
solutions to (\ref{bernoulli}) of cylindrical form, $a(z,t) = a_0(t)$,
one can do the integral to find
\begin{equation}
\ddot{a}_0\ln\left(\frac{4L^2}{a_0}\right) = \frac{\dot{a}_0^2}{2a_0}.
\label{cyl}
\end{equation}
It is easy to show that an an asymptotic solution of (\ref{cyl})
is given by
\begin{equation}
a_0 \propto \frac{\Delta t}{\ln(\Delta t)^{1/2}},
\label{cav0}
\end{equation}
corresponding to a power law with a small logarithmic correction.
Indeed, initial theories of bubble pinch-off \cite{LKL91,OP93}
treated the case of an approximately cylindrical cavity,
which leads to the radial exponent $\alpha=1/2$, with
logarithmic corrections.

\begin{figure}
\centering
\includegraphics[width=8cm]{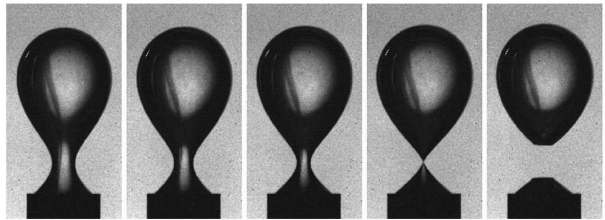}
\caption{
The pinch-off of an air bubble in water \cite{TET07}. An
initially smooth shape develops a localised pinch-point.
\label{shape}   }
\end{figure}

However both experiment \cite{TET07} and simulation \cite{EFLS07}
show that the cylindrical solution is unstable; rather, the pinch
region is rather localised, see Fig.~\ref{shape}. Therefore, it is
not enough to treat the width of the cavity as a constant $L$; the
width $\Delta$ is itself a time-dependent quantity. In \cite{EFLS07}
we show that to leading order the time evolution of the integral
equation (\ref{bernoulli}) can be reduced to a set of ordinary
differential equations for the minimum $a_0$ of $a(z,t)$, as well as
its curvature $a_0''$.

\begin{figure}
\centering
\psfrag{t'}{\Large$t'$}
\includegraphics[width=7cm]{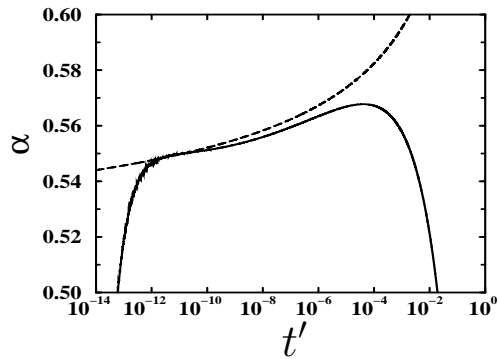}
\caption{A comparison of the exponent $\alpha$ between full
numerical simulations of bubble pinch-off (solid line) and the
leading order asymptotic theory \eref{asymp-exp} (dashed line).
\label{compare}   }
\end{figure}

Namely, the integral in (\ref{bernoulli}) is dominated by a
local contribution from the pinch region. To estimate this
contribution, it is sufficient to expand the profile around the
minimum at $z=0$: $a(z,t) = a_0 + a''_0/2 z^2 + O(z^4)$. As
in previous theories, the integral depends logarithmically on
$a$, but the axial length scale is provided by the inverse
curvature $\Delta\equiv (2a_0/a''_0)^{1/2}$. Thus evaluating
(\ref{bernoulli}) at the minimum, one obtains \cite{EFLS07}
to leading order
\begin{equation}
\ddot{a}_0\ln(4\Delta^2/a_0) = \dot{a}_0^2/(2a_0),
\label{a0}
\end{equation}
which is a coupled equation for $a_0$ and $\Delta$.
Thus, a second equation is needed to close the system,
which is obtained by evaluating the the second derivative of
(\ref{bernoulli}) at the pinch point:
\begin{equation}
\ddot{a}''_0\ln\left(\frac{8}{e^3a''_0}\right)
-2\frac{\ddot{a}_0a''_0}{a_0} = \frac{\dot{a}_0\dot{a}_0''}{a_0} -
\frac{\dot{a}_0^2a_0''}{2a_0^2}.
\label{Deltaequ}
\end{equation}

The two coupled equations (\ref{a0}),(\ref{Deltaequ})
are most easily recast in terms of the
time-dependent exponents
\begin{equation}
2\alpha\equiv -\partial_{\tau}a_0/a_0, \quad
2\delta \equiv -\partial_{\tau}a''_0/a''_0,
\label{exp}
\end{equation}
where $\tau \equiv -\ln t'$ and $\beta = \alpha-\delta$,
are a generalisation of the usual exponents $\alpha$ and
$\beta$. The exponent $\delta$ characterises the time dependence
of the aspect ratio $\Delta$. Returning to the collapse (\ref{cyl})
predicted for a constant solution, one finds that $\alpha = 1/2$
and $\delta = 0$. In the spirit of the the previous subsection,
this is the fixed point corresponding to the cylindrical solution.
Now we expand the values of
$\alpha$ and $\delta$ around their expected asymptotic values
$1/2$ and $0$:
\begin{equation}
\alpha=1/2+u(\tau),\quad
\delta =v(\tau).
\label{exp_ad}
\end{equation}

To leading order in $\delta$, the resulting equations are
\begin{equation}
\partial_{\tau}u = -8vu^2, \quad \partial_{\tau}v = -8v^3,
\label{deg_ad}
\end{equation}
which describe perturbations around the leading-order similarity
solution. These equations are
analogous to \eref{a1}, but they have a degeneracy of third
order, rather than second order. Equations \eref{deg_ad} are
easily solved to yield, in an expansion for small $\delta$
\cite{EFLS07},
\begin{equation}
\alpha=1/2+\frac{1}{4\sqrt{\tau}}+O(\tau), \quad
\delta =\frac{1}{4\sqrt{\tau}} + O(\tau^{-3/2}).
\label{asymp-exp}
\end{equation}
Thus the exponents converge toward their asymptotic values
$\alpha=\beta=1/2$ only very slowly, as illustrated in Fig.~\ref{compare}.
This explains why typical experimental values are found in the
range $\alpha\approx 0.54-0.58$ \cite{TET07}, and why there is
a weak dependence on initial conditions \cite{BMSSPL06}.

The cubic equation (\ref{exp_ad}) applies to the exponents
$\alpha$,$\delta$, rather than the solution itself, as in the
previous subsection, where we dealt with mean curvature flow.
In fact, it is easy to re-analyse the solution to the mean curvature
problem, and to formulate it in terms of time-dependent exponents.
Let us define $a_0$ and $a_0''$ for the mean curvature problem
through $h^2(z,t) = a_0 + a''_0/2 z^2 + O(z^4)$. From (\ref{blmc})
it follows that $a_0 = 2t'$ and $a_0''=2(\sqrt{2}+g)g''$.
This gives $a_0'' = 4\sqrt{2} a_1$, and thus
\begin{equation}
\partial_{\tau}\delta = -2\delta^2,
\label{deg_mc}
\end{equation}
instead of (\ref{deg_ad}). This means mean curvature flow,
formulated in terms of exponents, once more gives a
quadratic non-linearity, rather than the cubic term (\ref{deg_ad})
found for bubble collapse. From (\ref{deg_mc})
one finds to leading order $\delta_{mc} = 1/(2\tau)$ and
$\alpha_{mc} = 1/2 + O(1/\tau^2)$ for the mean curvature flow.

\section{Limit cycles}
\label{limit}
An example for this kind of blow-up was introduced into the
literature in \cite{C93} in the context of cosmology. There is
considerable numerical evidence \cite{Gund03} that discrete
self-similarity occurs at the mass threshold for the formation of
a black hole. The same type of self-similarity
has also been proposed for singularities of the Euler equation \cite{PS05}
and for a variety of other phenomena \cite{Sor98}.
A reformulation of the original cosmological problem leads to
the following system:

\begin{eqnarray} 
&& f_x=\frac{(a^2-1)f}{x}, \label{cosm:a}  \\
&& (a^{-2})_x=\frac{1-(1+U^2+V^2)/a^2}{x}, \label{cosm:b}  \\
&& (a^{-2})_t=\left[\frac{(f+x)U^2-(f-x)V^2}{x}+1\right]/a^2-1,
\label{cosm:c}\\
&& U_x=\frac{f[(1-a^2)U+V]-xU_t}{x(f+x)}, \label{cosm:d} \\
&& V_x=\frac{f[(1-a^2)U+V]+xV_t}{x(f-x)}. \label{cosm:e}
\end{eqnarray}

In \cite{MG03}, the self-similar description corresponding
to the system (\ref{cosm:a})-(\ref{cosm:e}) was solved using formal
asymptotics and numerical shooting procedures. This leads to the
solutions observed in \cite{C93}.
Below we propose a very similar system, which we solve analytically,
and which shows the same type of limit cycle behaviour:

\begin{eqnarray}
&& u_t + u_x = 2fv, \label{f:a}  \\
&& v_t + v_x = -2fu, \label{f:b}  \\
&& f_t = f^2. \label{f:c}
\end{eqnarray}

The simplified system (\ref{f:a})-(\ref{f:c}) 
can be solved introducing characteristics:
\begin{equation}
\eta=x+t, \quad \nu=x-t,
\label{char1}
\end{equation}
which leads to
\begin{equation}
u_{\eta}=fv, \quad v_{\eta}=-fu.
\label{char2}
\end{equation}
The system (\ref{char2}) has to be solved at constant $\nu$ with
initial conditions
\begin{equation}
u(\nu,\nu)=u_{init}(\nu), \quad v(\nu,\nu)=v_{init}(\nu),
\label{c_init}
\end{equation}
where $u_{init}(x)\equiv u(x,0)$ and $v_{init}(x)\equiv v(x,0)$
are the profiles at $t=0$.

The solution of (\ref{f:c}) is
\begin{equation}
f(x,t)=\frac{1}{1/f_{init}(x)-t},
\label{fsol}
\end{equation}
where $f_{init}(x)$ is the initial profile, which
can be written in terms of $\nu,\eta$ as
\begin{equation}
f(\eta,\nu) = \frac{1}{1/f_{init}((\eta+\nu)/2)-(\eta-\nu)/2}.
\label{fsol1}
\end{equation}
The function $f$ can be eliminated from (\ref{char2}) using
the transformation
\begin{equation}
\zeta = \int_{\nu}^{\eta} f(\eta',\nu) d\eta',
\label{v_trans}
\end{equation}
which results in the system for a harmonic oscillator:
\begin{equation}
u_{\zeta}=v, \quad v_{\zeta}=-u.
\label{char3}
\end{equation}
The solution is
\begin{equation}
u = C(\nu)\sin(\zeta+\phi_0(\nu)),\quad
v = C(\nu)\cos(\zeta+\phi_0(\nu)),
\label{sin_cos}
\end{equation}
with $C(\nu) = \sqrt{u_{init}^2(\nu)+v_{init}^2(\nu)}$ and
$\phi_0(\nu) = \arcsin(u_{init}(\nu)/C(\nu))$.

According to (\ref{fsol}), a singularity of the PDE system 
(\ref{f:a})-(\ref{f:c}) 
first occurs at the {\it maximum} $f_0$ of $f_{init}(x)$, which
we assume to occur at $x=0$ without loss of generality.
Thus locally we can write $f(x,0) \approx f_0-ax^2$, and for small
$t'=t_0-t$ we have
\begin{equation}
f=\frac{t'^{-1}}{1+at_0^2\xi^2} , \quad \xi=\frac{x}{t'^{1/2}},
\label{fs}
\end{equation}
for any finite $\xi$.
Using this explicit form of $f$, (\ref{v_trans}) can be integrated
to find $\zeta$. At constant $\xi$ we have
\begin{equation}
\nu = -t_0+\xi t'^{1/2}+t', \quad \eta = t_0+\xi t'^{1/2}-t',
\label{ne_trans}
\end{equation}
so taking the limit $t'\rightarrow 0$, the leading order result for
(\ref{v_trans}) is
\begin{equation}
\zeta = -\ln\left[\frac{t'(1+at_0^2\xi^2)(1+at_0^2)}{t_0}\right]
\equiv \tau + \phi(\xi).
\label{zeta_res}
\end{equation}
Thus as the singularity is reached, $t'\rightarrow 0$, the variable
$\zeta$ goes to infinity. This means the singularity corresponds to
the long-time limit of the dynamical system (\ref{char3}), which
will perform a harmonic motion according to (\ref{sin_cos}).
Namely, for $t'\rightarrow 0$ the result is
\begin{eqnarray}
&& u = \sqrt{u_{init}^2(-t_0)+v_{init}^2(-t_0)}\sin\left\{
\arcsin\left(\frac{u_{init}(-t_0)}
{\sqrt{u_{init}^2(-t_0)+v_{init}^2(-t_0)}}\right)
\right. \nonumber \\
&& \left. - \ln\left[\frac{t'(1+at_0^2\xi^2)(1+at_0^2)}{t_0}\right]\right\}.
\label{ses}
\end{eqnarray}

Thus the singular solution is of the general form
\begin{equation}
u = \psi(\phi (\xi )+\tau ),
\label{dss_structure}
\end{equation}
where $\psi$ is periodic in $\tau$. This is a particularly
simple version of discretely self-similar behaviour. Note
that the character of the solution is different from the
travelling waves described in section \ref{travel}.

\section{Strange attractors and exotic behaviour}
In connection to limit cycles and in the context of singularities
in relativity, a few interesting situations have been found numerically
quite recently. One of them is the existence of Hopf bifurcations where a
self-similar solution (a stable fixed point) is transformed into a discrete
self-similar solution (limit cycle) as a certain parameter varies (see
\cite{HE97}). Other kinds of bifurcations, for example of the
Shilnikov type, are found as well \cite{ABT06}. Now we demonstrate
that chaotic behaviour is also possible.

In section \ref{quad} we treated a system of an infinite number
of ordinary differential equations for the coefficients of the
expansion of an arbitrary perturbation to an explicit solution.
Such high-dimensional systems in principle allow for a rich
variety of dynamical behaviours, including those
found in classical finite dimensional dynamical systems, such as
chaos.
Consider for instance an equation for the perturbation $g$ (the
analogue of (\ref{gequ})) of the form
\begin{equation}
g_{\tau }=\mathit{L}g+F(g,g),
\label{ggen}
\end{equation}
where  $\mathit{L}g$ is a linear operator. Assuming an
appropriate non-linear structure for the function $F$, an
arbitrary nonlinear (chaotic) dynamics can be added.

To give an explicit example of a system of PDEs exhibiting
chaotic dynamics, consider the structure of the example given
in the previous section \ref{limit}. It can be generalised to produce {\it any}
low-dimensional dynamics near the singularity. Namely, let us
generalise the system (\ref{f:a})-(\ref{f:c}) to

\begin{eqnarray}
&& u^{(i)}_t + u^{(i)}_x = 2fF_i(\{u^{(i)}\}), \quad i=1,\dots,n,
\label{ch:a}  \\
&& f_t = f^2. \label{ch:d}
\end{eqnarray}

Using the transformation (\ref{char1}), the first $n$ equations
are turned into:
\begin{equation}
u^{(i)}_{\eta}/f=F_i, \quad i=1,\dots,n,
\label{sys1}
\end{equation}
which is an ODE system for constant $\nu$. The system (\ref{sys1})
has to be solved with initial conditions
\begin{equation}
u^{(i)}(\nu,\nu)=u^{(i)}_{init}(\nu), \quad \quad i=1,\dots,n.
\label{sys_init}
\end{equation}

As before, the function $f$ can be eliminated using (\ref{v_trans}),
resulting in the general non-linear dynamical system
\begin{equation}
u^{(i)}_{\zeta}=F_i, \quad i=1,\dots,n.
\label{sys2}
\end{equation}
Now if one chooses $n=3$ and
\begin{equation}
F_1 = \sigma(u^{(2)}-u^{(1)}), \quad
F_2 = \rho u^{(1)}-u^{(2)}-u^{(1)}u^{(3)}, \quad
F_3 = u^{(1)}u^{(2)}-\beta u^{(3)},
\label{sys_L}
\end{equation}
(\ref{sys2}) is the Lorenz system \cite{Strogatz94}.

As before, for $t'\rightarrow 0$, the variable $\zeta$ goes to
infinity, and near the singularity one is exploring the long-time
behaviour of the dynamical system (\ref{sys2}). In the case of
(\ref{sys_L}), and for sufficiently large $\rho$, the resulting
dynamics will be chaotic. Specifically, taking $\sigma=10$,
$\rho=28$, and $\beta=8/3$, as done by Lorenz \cite{L63}, the
maximal Lyapunov exponent is $0.906$.
Now if $U^{(1)}(\nu,\zeta),U^{(2)}(\nu,\zeta),U^{(3)}(\nu,\zeta)$
is a solution of (\ref{sys2}) with initial
conditions (\ref{sys_init}), the final form of the similarity
solution is
\begin{equation}
u^{(i)}(x,t) = U^{(i)}(-t_0+\xi t'^{1/2},\tau+\phi(\xi)), \quad i=1,2,3.
\label{chao_sim}
\end{equation}
In the limit $t'\rightarrow 0$ one cannot replace
the first argument by $-t_0+\xi t'^{1/2}\approx -t_0$. The reason
is that in the limit of large $\tau$ two trajectories diverge like
$\exp(\lambda_{max}\tau)=t'^{-\lambda_{max}}$, where the largest
Lyapunov exponent is larger than $1/2$.

\section{Outlook}
The singularities described in this paper are point-like in the sense that
they occur at a single point $x_0$ at a given time $t_{0}$. There are two
important situations we have not discussed and for which the dynamical
systems point of view and analytical approach does not apply. First, the
case of singularities not developing at single points, but on sets of
finite measure. This is the case in a few simple examples of
reaction-diffusion equations of the family
\begin{equation}
u_{t}-\Delta u=u^{p}-b\left\vert \nabla u\right\vert^{q}\quad \mbox{for}
\quad x \in\Omega
\label{react-diff}
\end{equation}
where depending on values of $p>1$ and $q>1$ singularities in the form of
blowing-up $u$ may be regional ($u$ blows up in subsets of $\Omega $ of
finite measure) or even global (the solution blows-up in the whole domain).
See for instance \cite{Souplet} and references therein. Another interesting
possibility is that the stable fixed point which is approached depends on the
initial conditions. Thus exponents could vary either discretely or
continuously with the choice of initial conditions. The infinite sequence
of similarity solutions found for (\ref{visc}) is {\it not} an example
for such behaviour, since only one solution is stable. All other fixed
points will not be visible in practise.

Singularities may
even happen in sets of fractional Hausdorff dimension, i.e., fractals. This
is the case of the inviscid one-dimensional system for jet breakup
(cf. \cite{FV00}) and might be case of Navier-Stokes system in
three dimensions, where
the dimension of the singular set at the time of first blow-up is at most $1$
(cf. \cite{CKN82}). This connects with the second issue we did not address
here. It is the nature of the singular sets both in space and time. \ In
many instances, existence of global in time (for all $0\leq t<\infty $)
solutions to nonlinear problems can be established in a \textit{weak sense},
that is, allowing certain kind of singularities to develop both in space and
time. In the case of 3-D Navier-Stokes system, the impossibility of
singularities "moving" in time, that is of curves $\mathbf{x}=\mathbf{%
\varphi }(t)$ $\ $in the singular set is well-known \cite{CKN82}. Hence,
provided a certain singularity does not persist in time, the question
is how to continue the solutions after a singularity has developed.

\ack
This paper is an outgrowth of discussions between the authors and
R. Deegan, preparing a workshop on singularities at the
Isaac Newton Institute, Cambridge. We thank J. M. Martin-Garcia and J. J.
L. Velazquez for fruitful discussions and for providing us with valuable
references.

\section*{References}
\providecommand{\newblock}{}

\end{document}